\documentclass{elsarticle}

\usepackage{amsmath}
\usepackage{amssymb}
\usepackage{graphicx}
\usepackage{hyperref}
\usepackage{color}
\usepackage{subcaption}

\newcommand{\umax}{u_{\text{max}}}
\newcommand{\be}{\begin{equation}}
\newcommand{\ee}{\end{equation}}
\newcommand{\partialderiv}[2]{\frac{\partial#1}{\partial#2}}
\newcommand{\diag}{\mathrm{diag}\,}

\begin{document}

\title{Optimal strategies for the control of autonomous vehicles in data
assimilation}

\author[ut]{D.~McDougall\corref{cor1}}
\ead{damon@ices.utexas.edu}

\author[njit]{R.O.~Moore}
\ead{rmoore@njit.edu}

\cortext[cor1]{Corresponding author}

\address[ut]{201 E. 24th St., Stop C0200\\
Institute for Computational Engineering and Sciences\\
University of Texas at Austin\\
Austin, TX, 78712}

\address[njit]{Department of Mathematical Sciences\\
New Jersey Institute of Technology\\
University Heights, Newark, NJ\, 07102}

\begin{abstract}
We propose a method to compute optimal control paths for autonomous vehicles deployed for the purpose of inferring a velocity field.  In addition to being advected by the flow, the vehicles are able to effect a fixed relative speed with arbitrary control over direction.  It is this direction that is used as the basis for the locally optimal control algorithm presented here, with objective formed from the variance trace of the expected posterior distribution.  We present results for linear flows near hyperbolic fixed points.
\end{abstract}

\begin{keyword}
  Bayesian inverse problem \sep Lagrangian data assimilation \sep Optimal control \sep Ocean glider

  \MSC[2010] 49M \sep 62F \sep 62L \sep 93C \sep 65C
\end{keyword}

\maketitle

\section{Introduction}

The need for a more accurate and better resolved estimate of oceanic flows is being driven by a number of pressing global issues, including the crisis facing many species of fish and waterborne organisms, the mitigation of pollutants resulting from spills and offshore contamination, and the important role played by ocean dynamics on climate change.
Scientific efforts to estimate ocean flow began in the 1980s with the work of Robinson~\cite{Robinson1981}, but has enjoyed limited success due to a lack of observational data.  In an effort to improve the current state of understanding of the world's oceans, autonomous vehicles (AVs) are being deployed for the collection of physical oceanography data in a growing number of projects around the globe.  One example of AVs are autonomous underwater vehicles (AUVs), which are equipped with adjustable wings that convert vertical momentum induced by battery-powered changes in buoyancy into forward momentum.

Notwithstanding a myriad of challenges caused by their immersion in a high-inertia environment and the limited power and bandwidth available for communication, the effectiveness of AVs comes from their ability to remain in the field for prolonged deployment and from their capacity for controlled self-propulsion.  To understand why this control is so important, one need simply consult the growing literature on filtering of complex systems.  The dimension associated with a meaningful estimate of an oceanic velocity field is far larger than a typical AV cohort size, so that the observations are necessarily sparse in space.  They may also be sparse in time for autonomous underwater vehicles (AUVs) that are unable to communicate while submerged.  Even large cohorts of Lagrangian drifters that simply advect with the flow are often subject to collapse onto a small subset of stable dynamic features of the flow, rendering their observations increasingly less informative with time.  Recent demonstrations that such drifters face inherent information barriers suggest that simply increasing the size of a deployed host is not effective without also implementing a method to allow them to cross flow barriers~\cite{chen_information_2014}.

To better understand this issue, we reconstruct a two-dimensional flow field from observations taken by AVs that are capable of limited locomotion in addition to advecting with the flow.  We refer to these AVs as {\em gliders}, in distinction to drifters. They are capable of measuring the surrounding fluid velocity either directly using Doppler velocimetry or indirectly using sequential position measurements, depending on the instrumentation on board.  In this work we focus on the former case and use the velocity observations to perform an ``identical twin'' experiment based on a perfect model, with the goal of constructing an estimate of the true
flow field function from observations, along with an associated uncertainty.
Our approach to control follows the underlying philosophy of \cite{mcdougall_jones_2015} that mesoscale
ocean flow fields are dominated by coherent features, such as jets and eddies, which not only dictate the sites where the most informative observations can be taken to estimate these structures, but also provide the most effective transport mechanisms for navigation of the physical domain.

If the ocean drifter can be controlled to move into and through these
structures then the information gained should be richer in terms of capturing
the key properties of the flow field.  Whereas in~\cite{mcdougall_jones_2015} this was shown by utilizing an
ad hoc control strategy, the present work attempts to put this philosophy on a more systematic footing
through the application of an \textit{optimal} control strategy.  Each optimal control calculation
is inserted between successive observations in a sequential filter, utilizing the posterior field estimate to
computing a control that most effectively minimizes the next expected posterior variance.



In general, reconstructing the flow requires one to solve an inverse
problem which is most naturally posed in a Bayesian framework~(\cite{Kalnay2002}), the solution to which is a probability distribution
on the appropriate function space or parameter space, depending on the velocity field's representation.  Kalman filters and Kalman smoothers~(\cite{Kalman1960,Kalman,Sorenson1960,Evensen2006,Houtekamer1998,
Anderson2003}) completely quantify this distribution by its first two moments in the case of linear processes and Gaussian initial distributions, and are popular approximate methods for nonlinear flows and non-Gaussian distributions.  The Kalman filter is the assimilation mechanism used here, ostensibly due to the use of linear stationary velocity fields and Gaussian distributions but also to provide the straightforward extension to weakly nonlinear flow models through the extended Kalman filter~\cite{Kalnay2002}.

Alternative approaches include variational methods that view the solution not as a distribution but as
the argument of a cost
function that is optimized~\cite{Lorenc2000,Bengtsson1975,Lewis1985,
Lorenc1986,LeDimet1986,Talagrand1987,Courtier1994,Lawless2005,Lawless2005a},
and particle filtering methods
(\cite{Doucet2001,VanLeeuwen2010}) that approximate the continuous posterior random
variable by a discrete set of state realizations (particles) with associated
weights.  Updating the particles and weights as new observations are made is difficult and
can lead to degenerate posteriors in high-dimensional problems~\cite{Bickel:2008bx}.  Sampling
methods utilize the Metropolis-Hastings algorithm
(\cite{Metropolis1953,Hastings1970}) to sequentially generate correlated samples
from the posterior distribution
(\cite{cotter2013mcmc,Cotter2009,Cotter2010,Lee2011,Apte2008a,Apte2007,Apte2008,
Herbei2009,Kaipio2000,Mosegaard1995,Roberts1997,Roberts1998,Roberts2001,
Beskos2009,Atchade2005,Atchade2006}).  Sampling methods are ideal for the online
computation of moments through the use of unbiased estimators, but are not
suited to applications where the data arrive sequentially, and are expensive in settings
where the distribution must be used in an inverse problem to determine the minimizing control.

Section~\ref{s:setup} describes the setup of the problem under consideration in a relatively broad setting, including the basic flow assumptions and assimilation model, followed by a derivation of the method used to find local minimizers of the objective.  Section~\ref{s:regularization} identifies an existence problem with the variational formulation and introduces a simple approach to address it.  Section~\ref{s:results} provides numerical results for four linear time-independent velocity fields modeling flow near hyperbolic fixed points with different stability properties.  Section~\ref{s:discussion} offers conclusions and discusses extensions to the methodology.



\section{Setup}
\label{s:setup}

The context of this study is oceanographic data assimilation, where the goal is to estimate a velocity field conditioned on sparse observational data recorded by Lagrangian {\em gliders\/} whose positions evolve according to the flow being assimilated in addition to a modest capacity for self-propulsion.
The gliders therefore move according to
\be
  \dot{z} = v(z(t), t) + u(t),
\ee
with $z\in \mathbb{R}^{dK}$ contains the positions of $K$ gliders in a $d$-dimensional domain (here, we take $d=2$), and the flow evolves according to a model expressed as
\be
\partialderiv{v}{t} = f(v).
\label{e:flowModel}
\ee
The glider's self-propulsion is expressed by control $u(t)$, where $u(t)\equiv 0$ represents the uncontrolled case (i.e., Lagrangian drifters~\cite{mariano_lagrangian_2002}) and $|u(t)|\equiv \umax$ represents the case considered here with the relative speed of gliders fixed at $\umax$ \cite{Leonard2010}.

The direct velocity observations considered here take the form
\be
y_i = v(z(t_i),t_i) + \eta_i,
\label{e:obs}
\ee
where $\eta_i\sim\mathcal{N}(0,R_i)$ is i.i.d Gaussian-distributed measurement noise.  Here, $y_i\in \mathbb{R}^{dK}$ and $v(z,t):=(v(z_1,t),v(z_2,t),\dots,v(z_K,t))^T$.  The mapping between ``state space'' (estimates of $v(x,t)$) and ``observation space'' (values of $y$) is referred to as the observation operator which, in this case, is formally given by the convolution
\be
H(v) = \delta(x-z(t_i))*v.
\label{e:obsOpVel}
\ee

If the observational data consist of glider positions, then an alternative to inferring velocities through, for example, dead reckoning before assimilation is to augment the state vector to include the glider positions~\cite{Kuznetsov2003}.  If the mean discretized velocity state estimate $\hat{v}$ is $dN$-dimensional, this produces an augmented state vector ${\tilde{v}}:=(\hat{z}^T,\hat{v}^T)^T$ of dimension $d(N+K)$, with the advantage that the observation operator is simply
\be
H(\tilde{v}) = \begin{pmatrix}I_{dK\times dK} & O_{dK\times NK}\end{pmatrix}\begin{pmatrix}\hat{z}\\ \hat{v}\end{pmatrix}.
\label{e:obsOpPos}
\ee

\subsection{Assimilation}

There exist a wide variety of techniques for the assimilation of Lagrangian observational data into a velocity model, depending on the context.  For linearly evolving flows, filters based on Gaussian or conditionally Gaussian~\cite{liptser_statistics_2001} processes such as the Kalman filter are appropriate, with suitable modification to reconcile Eulerian and Lagrangian descriptions~\cite{piterbarg_optimal_2008,Kuznetsov2003}.  These methods have been extended with some success to weakly nonlinear flows; fully nonlinear flows require alternative methods such as particle filters or variational frameworks to explore the posterior distribution.

The ultimate goal is to combine both the flow model in Eqn.~\ref{e:flowModel}
and the observations in Eqn.~\ref{e:obs} to obtain an estimate (with uncertainty) of the underlying
flow $v$.  In terms of the Bayesian formulation, we seek the posterior probability distribution of the flow field conditioned on noisy observations, i.e.,
\begin{equation}
  \frac{\mu^y}{\mu_0}(v) \propto \mathbb{P}(y | v),
\end{equation}
where the posterior distribution $\mathbb{P}(v | y)$ has measure $\mu^y$ and
the prior distribution $\mathbb{P}(v)$ has measure $\mu_0$.  This formulation
of Bayes's rule, in terms of the Radon-Nikodym derivative, is necessary to
ensure well-posedness of the posterior distribution on function
space~\cite{stuart2010inverse}.
In the case of a linear evolution model $f$, an initially Gaussian velocity field estimate maintains this property, allowing its complete characterization through two moments, the mean $\hat{v}$ and covariance $\hat{P}$, with
\begin{align}
  \hat{v}_i  &= \mathbb{E}(v | y_i) \\
  \hat{P}_i(x, x') &= \mathbb{E}( (v(x ) - \mathbb{E}(v(x )))
                                (v(x') - \mathbb{E}(v(x')))^\top).
\end{align}
Here, $\hat{v}_i$ and $\hat{P}_i$ incorporate observations up to time $t_i$.
They are updated sequentially upon the arrival of new data.  If observations
take the form expressed in Eqn.~\ref{e:obsOpVel} the update formula is given by
the Kalman filter equations~(\cite{ide_uni_1997}),
\begin{align}
  \hat{v}_i &= \hat{v}_{i-1} + K_i(y_i - \hat{v}_{i-1}(z(t_i))), \\
  \hat{P}_i &= \hat{P}_{i-1} + K_i H(v;u) \hat{P}_{i-1}.
\end{align}
Here $K_i$ is called the Kalman gain and is given by,
\begin{equation}
  K_i = \hat{P}_{i-1} H_i^\top ( H_i \hat{P}_{i-1} H_i^\top + R)^{-1}.
\end{equation}

\subsection{Controlling the drifters}

Our objective is to control gliders to where their measurements of the flow
will minimize the trace of the posterior covariance,
\begin{equation}
  \Gamma(z(t_i)) = \text{tr} \left( \hat{P}_{i-1} - \hat{P}_{i-1} H_i^\top
  (H_i \hat{P}_{i-1} H_i^\top + R)^{-1} H_i \hat{P}_{i-1} \right).
\end{equation}
The trace of the posterior covariance matrix can be thought of as the `global'
uncertainty.  Therefore minimising this quantity can be thought of as
maximising our knowledge of the underlying ocean flow, which is crucial for use
in a prediction.

To impose the fixed relative speed constraint, we take a similar approach to
Zermelo's problem~\cite{BrysonJr.1975} and replace the $2K$-dimensional control
$u$ with
\begin{equation}
u = \umax\hat{u} := \umax(\cos\theta_1, \sin\theta_1, \cos\theta_2, \dots,
\sin\theta_K)^{\top},
\end{equation}
reducing the number of degrees of freedom to $K$ angles $\theta_k$.  The
objective function we aim to minimize is therefore
\begin{equation}
  \Phi(z(t_j)) = -\Gamma(z(t_j)) -
  \int_{t_{j-1}}^{t_j}\lambda(t)^{\top}(\dot{z}-v(z)-\umax\hat{u}(t))\,dt.
\end{equation}
Setting the first variation of $\Phi$ with respect to $z(t)$ equal to zero gives the Euler-Lagrange equations,
\begin{align}
  \dot{z} &= v(z) + \umax \hat{u}(t), \quad z(t_{j-1}) = z_{j-1}, \\
  \dot{\lambda} &= -\left( \partialderiv{v}{z} \right)^{\top} \lambda, \quad
  \lambda(t_j) = -\left( \left. \partialderiv{\Gamma} z \right|_{t_j} \right)^{\top}.
  \label{e:costateEqn}
\end{align}
The costate $\lambda$ and control $\hat{u}$ are related through
\begin{equation}
  \lambda(t_j)^{\top} \partialderiv{\hat{u}}{\theta} = 0,
\end{equation}
where $\partial\hat{u}/\partial\theta$ is a block-diagonal matrix of 2-vectors
$(-\sin\theta_k, \cos\theta_k)^{\top}$, with size $2K\times K$.  Since each of these
2-vectors is orthogonal to the corresponding 2-vector component of $\hat{u}$,
this condition is equivalent to $2K$-vector $\lambda$ taking the form
\begin{equation}
  \lambda = -\begin{pmatrix}\tilde{\lambda}_1\hat{u}_1\\
  \tilde{\lambda}_2\hat{u}_2\\
  \vdots\end{pmatrix},
  \label{e:lambdaRep}
\end{equation}
essentially representing $\lambda$ as a vector of polar forms.

The above Euler-Lagrange system of state and costate equations naturally
suggests a shooting method to compute the minimizer of $\Phi$.  Here we take
an alternative approach and derive a relaxation method that is based on
the expression of the coupled one-way initial- and terminal-value problems for
the state and costate.  We arrive at a single second-order boundary value
problem that we reformulate as the steady state of a parabolic PDE.
Specifically, differentiating Eqn.~\ref{e:costateEqn} and using
Eqn.~\ref{e:lambdaRep}, we arrive at the equivalent equation for $K$-vector
$\theta$:
\begin{equation}
  \dot{\theta} = -\diag\hat{u}^{\top}\diag J^{\top}
  \left(\partialderiv{v}{z}\right)^{\top}\hat{u},
  \label{e:thetaEqn}
\end{equation}
where $\diag\hat{u}$ is the $2K\times K$ matrix with each 2-vector $u_k$ along the diagonal, and
\begin{equation}
  J = \begin{pmatrix}0 & -1 \\ 1 & 0\end{pmatrix}.
\end{equation}
Note that the $K$ costate amplitudes $\tilde{\lambda}_k$ have dropped out,
leaving just the angles $\theta_k$.  These angles are decoupled, i.e.,
Eqn.~\ref{e:thetaEqn} is just a $K$-fold copy of
\begin{equation}
  \dot{\theta}_k = -\begin{pmatrix}-\sin\theta_k & \cos\theta_k\end{pmatrix}
  \begin{pmatrix} \partialderiv{v_{kx}}{x_k} & \partialderiv{v_{ky}}{x_k} \\
  \partialderiv{v_{kx}}{y_k} & \partialderiv{v_{ky}}{y_k}\end{pmatrix}
  \begin{pmatrix}\cos\theta_k \\ \sin\theta_k\end{pmatrix}.
\end{equation}
The (forward) state equations and (backward) costate equations are
also decoupled, in addition to the boundary conditions for the state at
$t_{j-1}$, leaving as the only source of coupling the terminal condition for
$\lambda$ and therefore $\theta$. This strongly suggests that an iterative
numerical method using decoupled path optimization should be possible.  We will
proceed with the formulation of the relaxation method for a single AUV, with
the understanding that all AUV paths are coupled at $t_j$.

Differentiating Eqn.~\ref{e:costateEqn} and making use of Eqn.~\ref{e:thetaEqn},
we arrive at the boundary value problem
\begin{equation}
  \ddot{z} = \partialderiv{v}{z}\dot{z} - \frac1{\umax^2}
  J(\dot{z} - v)(\dot{z} - v)^{\top}J^{\top}
  \left(\partialderiv{v}{z}\right)^{\top}(\dot{z} - v).
  \label{eqn:2ndorder}
\end{equation}
The boundary condition at initial time $t_j$ is linear and Dirichlet:
\begin{equation}
  z(t_{j-1}) = z_{j-1}.
  \label{e:BCleft}
\end{equation}
At the terminal time $t_j$, however, the boundary condition is nonlinear unless $v$ is
linear and $\Gamma$ is quadratic:
\begin{equation}
  \dot{z}(t_j) = v(z(t_j)) + \umax\left.\frac{(\partial\Gamma/\partial z)^{\top}}
  {\|\partial\Gamma/\partial z\|}\right|_{z(t_j)}.
  \label{e:BCright}
\end{equation}
As noted earlier, this boundary condition also couples all glider positions at $t_j$ through the posterior covariance trace $\Gamma(z_1,\dots,z_K)$.

\subsubsection{Numerical method}

Equation~\ref{eqn:2ndorder} with boundary conditions~\ref{e:BCleft} and~\ref{e:BCright} expresses a two-point boundary value problem for which a number of solution techniques exist~\cite{keller_numerical_1990}.  We adopt a relaxation approach to finding solutions, in which we introduce an artificial time $\tau$ and formulate a partial differential equation (PDE) with solutions to the two-point boundary problem as its steady states:
\begin{equation}
  \partialderiv{z}{\tau} =
  \ddot{z} - \partialderiv{v}{z}\dot{z} + \frac1{\umax^2}
  J(\dot{z} - v)(\dot{z} - v)^{\top}J^{\top}
  \left(\partialderiv{v}{z}\right)^{\top}(\dot{z} - v).
  \label{e:relaxDiff}
\end{equation}
This PDE is solved with boundary conditions as above and with a simple initial condition such as the uniform state $z(t)\equiv z_{j-1}$.  The PDE is integrated in $\tau$ until the residual of Eqn.~\ref{eqn:2ndorder}, i.e., $\|\partial{z}/\partial \tau\|$, drops below a threshold value.  We use central differences in $t$ and a semi-implicit method in $\tau$, with $\ddot{z}$ computed implicitly and the remaining terms on the right-hand side of Eqn.~\ref{e:relaxDiff} computed explicitly.

The minimization problem posed above is only expected to have unique solutions under very specific conditions, such as linear velocity and quadratic objective, so the initial condition chosen for the relaxation method will play a large role in determining the local minimum found.  This is discussed in Sec.~\ref{s:discussion}.

As noted earlier, these equations are uncoupled apart from the boundary
condition at $t_j$.  We therefore adopt the following iterative method defined
on a uniform partition of $(t_{j-1},t_j]$:
\begin{enumerate}
  \item Initialize all paths $z_k(t,0)$, e.g., using constants $z_k(t,0)\equiv z_{k,j-1}$ or straight-line paths
    connecting $z_{k,j-1}$ to nearby maxima of $\Gamma$.
  \item \label{item1} For each AUV $k=1,\dots,K$, update the path by
    integrating Eqn.~\ref{e:relaxDiff} until the residual (i.e., right-hand
    side) has dropped below a prescribed threshold, using $\Gamma$ evaluated at
    fixed $z_l$, $l=1,\dots,k-1,k+1,\dots,K$ (using the most recently computed
    path iterate available).
  \item Repeat \ref{item1} until all updated residuals are below threshold.
\end{enumerate}

\section{Convexity and regularization}
\label{s:regularization}

It is clear that convexity for this optimization problem depends entirely on
the nature of the velocity field $v(z)$ and cost function $\Gamma(z(t_j))$.
Nontrivial $v$ and $\Gamma$ are expected to have multiple local minimizers and,
more critically, are not guaranteed to admit a solution within the class of
admissible trajectories $z(t)$.  The focus of this work is simply to find a
reasonable local minimizer as dictated by the choice of initial condition for
the relaxation method, so the discussion of minimizer selection is deferred to
future work.  It is possible to find scenarios even in the case of constant velocity and quadratic cost where a solution fails to exist, which our numerical studies corroborate; we analyse this scenario below to motivate a regularization of the terminal
condition~\ref{e:BCright}.

If $v\equiv v_0$ is constant, relaxation method~\ref{e:relaxDiff} collapses to the simple heat equation,
\be
z_\tau = \ddot{z},\quad z(t_0)=z_0,
\label{e:heateqn}
\ee
where assuming $\Gamma(z) = -\frac12z^TBz$ gives
\be
\dot{z}(t_1) =v_0 - \umax\frac{(B+B^T)z(t_1)}{\|(B+B^T)z(t_1)\|}.
\label{e:BCreduced}
\ee
Coupling between the gliders is quantified by off-diagonal elements, i.e.,
elements outside the K $2\times 2$ submatrices that lie along the
main diagonal.  To motivate the regularization below it is sufficient to
consider decoupled gliders with $B=I$, although the same analysis is possible
for any $B$ that is symmetric and positive-definite.  The boundary
condition~\ref{e:BCreduced} is then simply
\be
\dot{z}(t_1) =v_0 - \umax\frac{z(t_1)}{\|z(t_1)\|}
\ee
and steady states of Eqn.~\ref{e:heateqn} must satisfy
\be
z(t_1) = z(t_0) + (t_1-t_0)\left(v_0-\umax\frac{z(t_1)}{\|z(t_1)\|}\right).
\label{e:condition}
\ee
\begin{figure}[h!]
\begin{subfigure}{0.49\linewidth}
\begin{center}
\includegraphics[scale=0.35]{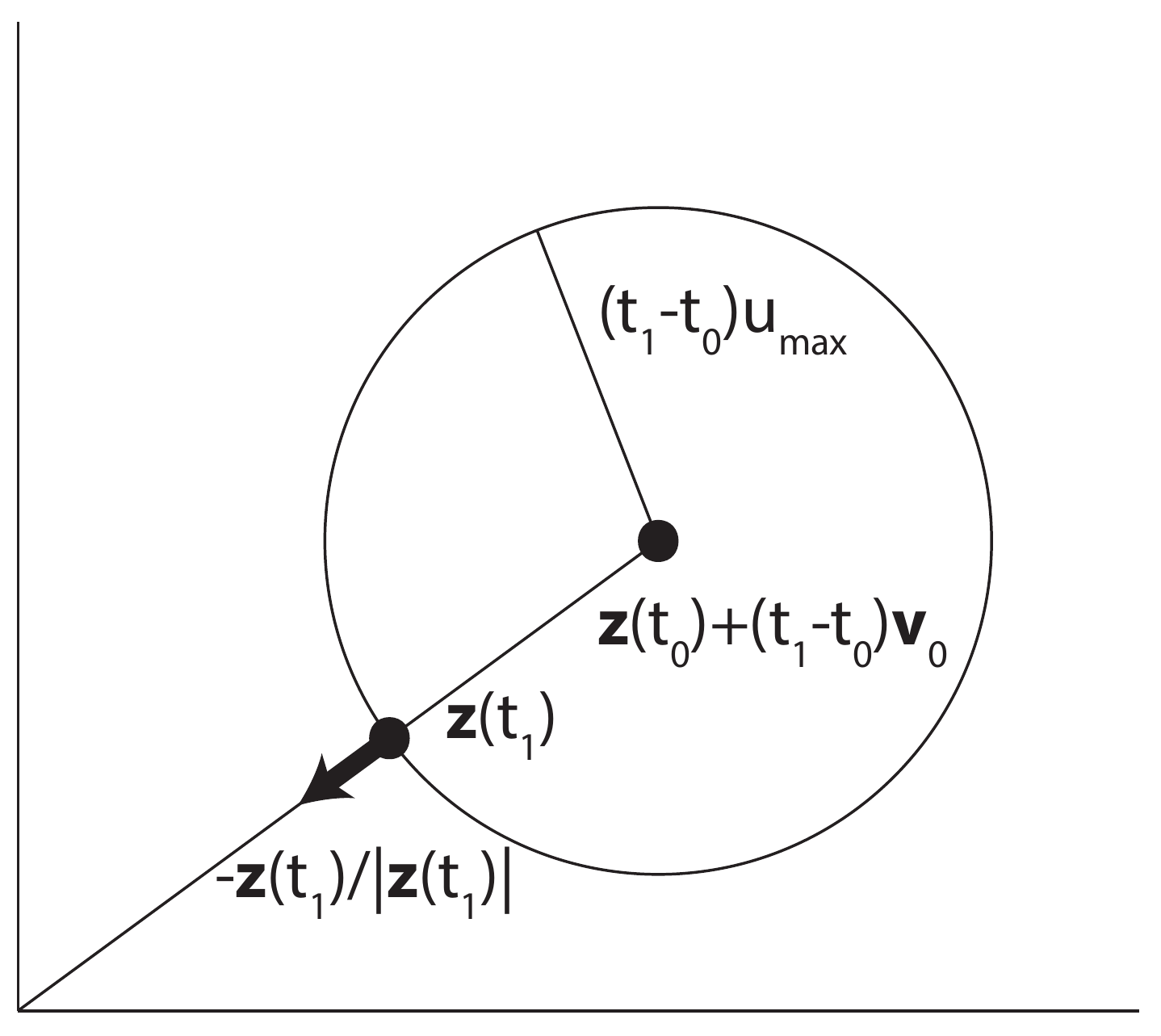}
\caption{\label{f:regularization:A}Unique solution.}
\end{center}
\end{subfigure}
\begin{subfigure}{0.49\linewidth}
\begin{center}
\includegraphics[scale=0.35]{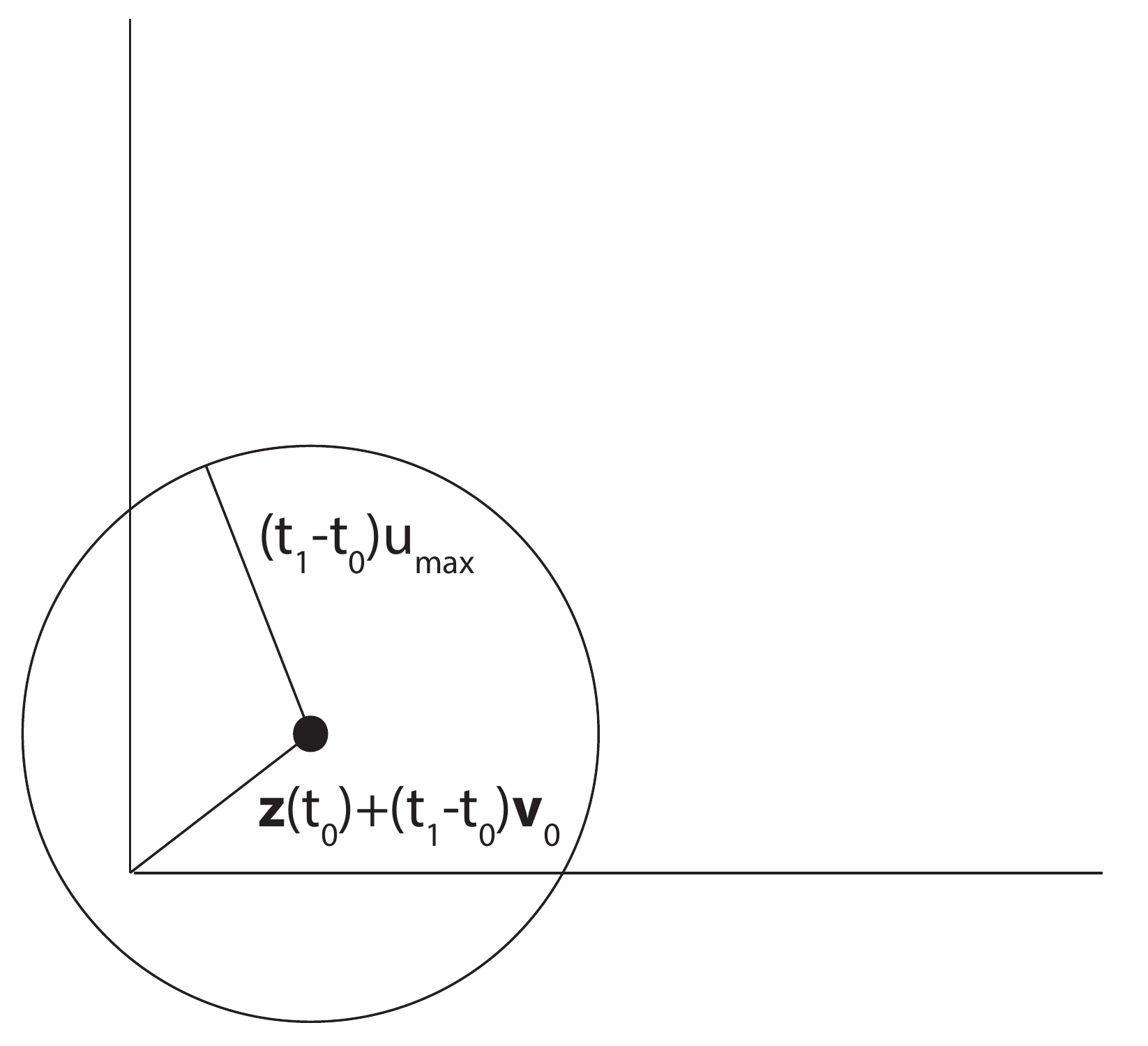}
\caption{\label{f:regularization:B}No solution.}
\end{center}
\end{subfigure}
\caption{\label{f:regularization}(A) Illustration of unique solution to condition~\ref{e:condition} when $\|z(t_0)+(t_1-t_0)v_0\|>(t_1-t_0)\umax$. (B) Illustration that no solution exists when $\|z(t_0)+(t_1-t_0)v_0\|<(t_1-t_0)\umax$.}
\end{figure}
Figure~\ref{f:regularization} illustrates that the existence and uniqueness of solutions in this case depends entirely on the difference between the uncontrolled distance from the origin (i.e., the global maximum of $\Gamma$) at time $t_1$ and the maximum relative distance possible under control $\umax$, i.e.,
\be
\delta L=\|z(t_0)+(t_1-t_0)v_0\| - (t_1-t_0)\umax.
\ee
If $\delta L > 0$, as in Fig.~\ref{f:regularization:A}, a solution exists and is unique.  If $\delta L < 0$, Fig.~\ref{f:regularization:B} demonstrates that no solution exists.  If $\delta L = 0$, the obvious solution has the glider land exactly at the maximum of $\Gamma$, although the boundary condition~\ref{e:BCreduced} is not well-defined here.

The failure to find a solution is directly related to the fixed-speed constraint on each glider; if its speed is large enough to send it past a nearby local maximum of the objective, this implies the lack of a minimizer in the space of smooth functions.  This problem can be addressed either by including piecewise smooth curves in the class of admissible paths or by replacing the fixed-speed constraint with a speed penalty in the objective function.  To minimize the impact on paths not affected by the problem discussed here, we adopt an intermediate approach of altering boundary condition~\ref{e:BCright} to allow violations of the fixed-speed constraint for gliders close to a local maximum, where the gradient of the objective drops below a specified magnitude.  The regularized boundary condition is as follows:
\be
  \dot{z}(t_j) = v(z(t_j)) + \begin{cases}
   \umax\left.\frac{(\partial\Gamma/\partial z)^{\top}}
  {\|\partial\Gamma/\partial z\|}\right|_{z(t_j)} & \mbox{ if } \|\partial\Gamma/\partial z\|>\gamma,\\
  \frac{\umax}{\gamma}\left.(\partial\Gamma/\partial z)^{\top}\right|_{z(t_j)} & \mbox{ otherwise.}
  \end{cases}
  \label{e:BCright_reg}
\ee

The impact of this regularization on the present case is that condition~\ref{e:condition} is replaced by
\be
z(t_1) = z(t_0) + (t_1-t_0)\left(v_0-\frac{\umax}{\gamma}z(t_1)\right),
\label{e:newcondition}
\ee
which has solution
\be
z(t_1) = \frac{\gamma}{\gamma + (t_1-t_0)\umax}\left(z(t_0)+(t_1-t_0)v_0\right).
\ee
This solution is only relevant if $\|\partial\Gamma/\partial z\|<\gamma$ at this value of $z$, i.e., if
$\gamma > \delta L$, which is always true in the case where no solution exists in the absence of regularization.
This solution limits to $z(t_1)=0$ (i.e., to the global maximum of $\Gamma$) as $\gamma\rightarrow 0$ and to $z(t_0)+(t_1-t_0)v_0$ (i.e., to the uncontrolled case) as $\gamma\rightarrow\infty$.

In practice, the value of $\gamma$ is chosen so to be small enough that it has no impact on the majority of trajectories that fall into the category of Fig.~\ref{f:regularization:A}, but large enough that it successfully resolves all cases suffering from a lack of existence due to the scenario depicted in Fig.~\ref{f:regularization:B}. 

\section{Results}
\label{s:results}


To test our method we consider a stationary linear flow
\[
v(z) = v_0 + Ax,\quad\mbox{with}\quad A \in {\mathbb R}^{2\times 2}.
\]
The state vector $(v_{0x}, v_{0y}, A_{11}, A_{12}, A_{21}, A_{22})^T$
is 6-dimensional with observation operator
\be
H_e = \begin{pmatrix} 1 & 0 & x_1 & y_1 & 0 & 0\\
0 & 1 & 0 & 0 & x_1 & y_1\\
& & \vdots & & & \\
1 & 0 & x_K & y_K & 0 & 0\\
0 & 1 & 0 & 0 & x_K & y_K
\end{pmatrix},
\label{e:obsOpLinear}
\ee
mapping from ${\mathbb R}^6$ to the set of $2K$ glider velocities.
Numerical simulations were performed to infer this linear velocity field using four different ``ground truths''.  In each case, the (unknown) constant component was given by $v_0=(1/2,\,-1/2)^T$, but the velocity's Jacobian $A$ was chosen to produce (a) a center, (b) an unstable node, (c) a saddle, and (d) a stable node.  Each case is considered below, where the numerical simulations were generated with cohorts of 1, 2, 5, and 10 gliders, each taking 100 sequential velocity observations with time interval of 0.1 units between successive observations.  The estimate was initialized at zero with variance of $10^6$.  Observational noise was set at unit variance.  The maximum glider speed was set to $\umax=1$.  In order to assess performance, all runs were performed with optimal control obtained using Eqn.~\ref{e:relaxDiff}, with no control, and with a fixed control for each assimilation stage with randomly selected direction.  Computations of the optimal control that failed to converge were redone with the regularization discussed in Sec.~\ref{s:regularization} using
\be
\gamma = \umax (t_j-t_{j-1})\|\frac{\partial^2\Gamma}{\partial z^2}\|.
\ee

\subsection{Flow about center}

\begin{figure}[h!]
\begin{subfigure}{0.49\linewidth}
\begin{center}
\includegraphics[scale=0.35]{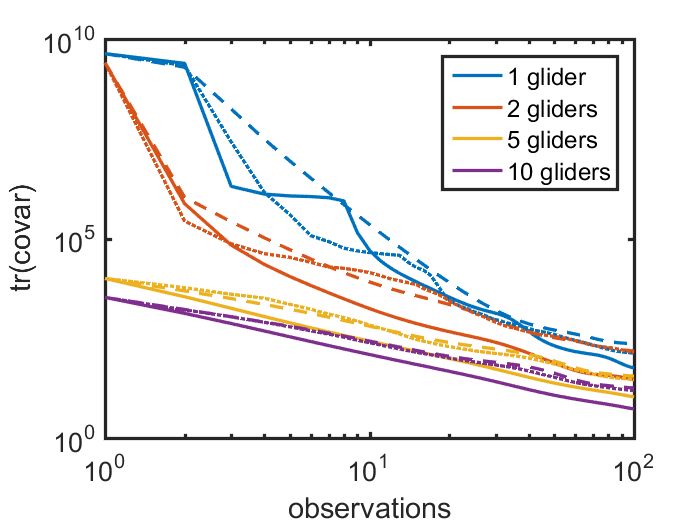}
\caption{\label{f:traceCenter:A}Covariance trace, flow about center.}
\end{center}
\end{subfigure}
\begin{subfigure}{0.49\linewidth}
\begin{center}
\includegraphics[scale=0.35]{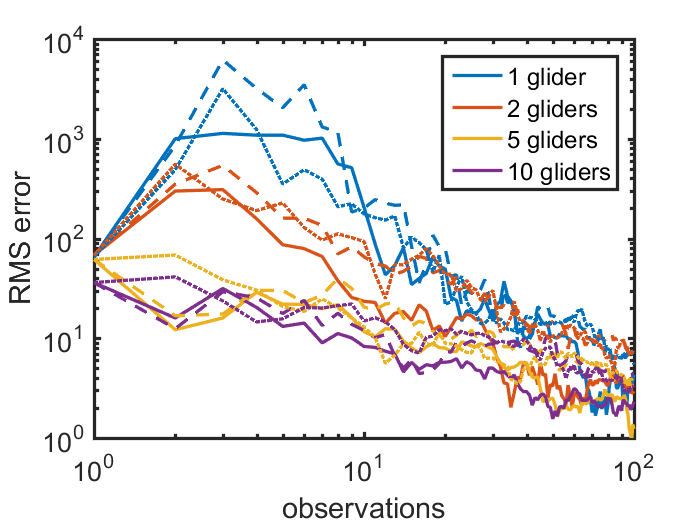}
\caption{\label{f:traceCenter:B}RMS error, flow about center.}
\end{center}
\end{subfigure}
\caption{\label{f:traceCenter}(A) Decrease in covariance trace with number of observations for flow about center.  For each glider cohort (blue, red, yellow, and purple curves respectively represent 1-, 2-, 5-, and 10-glider cohorts), the solid, dashed, and dotted curves respectively represent locally optimal control, no control, and random control.  Optimal control performs systematically better, particularly as the number of observations increases. (B) Same as in (A), but with RMS error used as the performance metric.  Initial growth in the RMS error for the 1- and 2-glider cohorts reflects the high impact of noisy observations on these underresolved cases.}
\end{figure}
\begin{figure}[h!]
\begin{subfigure}{0.49\linewidth}
\begin{center}
\includegraphics[scale=0.35]{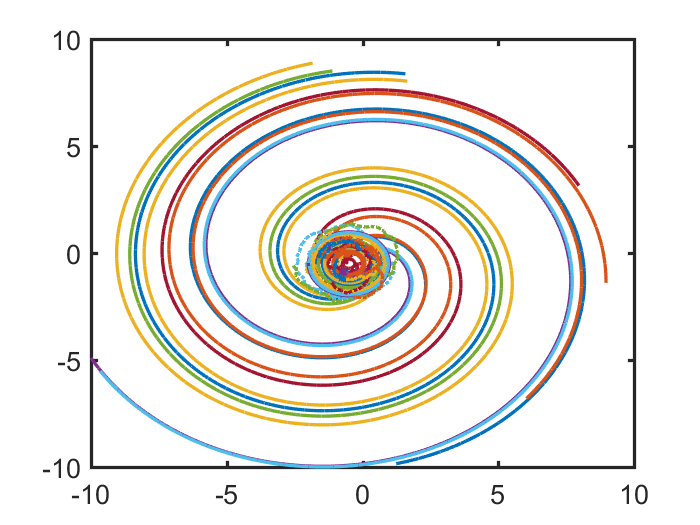}
\caption{\label{f:gliderPathsCenter:A}Optimal paths, flow about center.}
\end{center}
\end{subfigure}
\begin{subfigure}{0.49\linewidth}
\begin{center}
\includegraphics[scale=0.35]{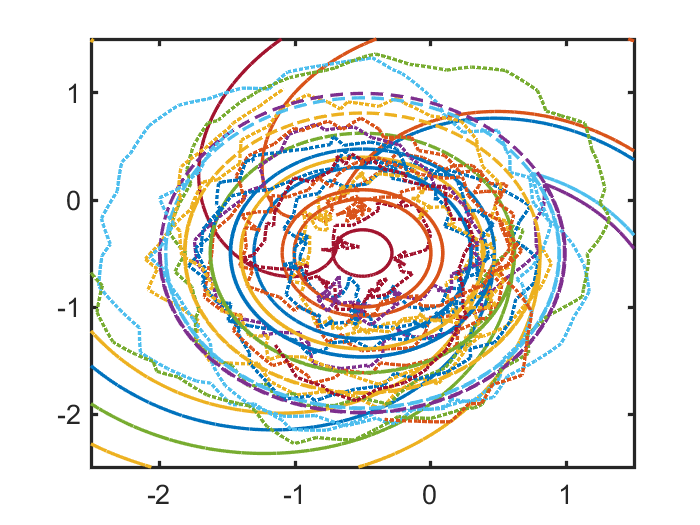}
\caption{\label{f:gliderPathsCenter:B}Blowup of (A).}
\end{center}
\end{subfigure}
\caption{\label{f:gliderPathsCenter}(A) Optimally controlled glider velocities for flow about center.  Solid lines are optimally controlled, dashed lines have no control, and dotted lines have random control. (B) Same as in (A) but with axes rescaled to give better picture of uncontrolled and randomly controlled paths.}
\end{figure}
The simulations performed around a center use a velocity field given by
\be
v_0 = \begin{pmatrix}1/2\\-1/2\end{pmatrix},\qquad A = \begin{pmatrix}0 & 1\\-1 & 0\end{pmatrix},
\ee
such that the center is located at $(-1/2,\,-1/2)$.
Figure~\ref{f:traceCenter} demonstrates the advantage of using optimal control, represented by solid lines, over no control, represented by dashed lines, for cohort sizes of 1, 2, 5, and 10 in terms of the covariance trace and the root-mean-square (RMS) error.  All cases show systematic improvement with the use of optimal control, with this improvement varying according to the usable information contained in the velocity field estimate.  The Kalman filter estimates have necessarily decreasing covariance trace while the RMS error shows an initial increase as the first few assimilations place too much confidence on random point estimates.
In the 1- and 2-glider cohorts, there is little systematic improvement in the estimate quality for the first few iterations since the estimate is unable to accurately inform the optimization computation.  Once a sufficiently reliable velocity estimate is established, however, the improvement is substantial.  This is particularly striking in the case of the single-glider cohort, where the RMS error drops precipitously between the 60th and 70th observation, breaking away from its uncontrolled counterpart.  Also included as dotted lines in Fig.~\ref{f:traceCenter} are the assimilation results using randomly controlled trajectories, where the direction of control is drawn from a uniform distribution of angles on $[0,2\pi)$ and held constant between observations.  The randomly controlled trajectories are seen to offer some benefit in the intermediate observation times where the velocity estimate is poor, but they ultimately perform almost identically to the uncontrolled trajectories, and considerably worse than the optimally controlled trajectories.

To help understand the performance improvement offered by the optimally controlled gliders, their trajectories are plotted in Fig.~\ref{f:gliderPathsCenter}~(A) and~(B).  As seen in Fig.~(B), the gliders simply orbit around the center in the absence of control, with similar albeit erratic trajectories in the case of random control.  The optimal control, however, pushes the gliders systematically toward larger distances from the origin.  This is due to the fact that the information provided by observation operator~\ref{e:obsOpLinear} increases with larger values of $x$ and $y$.  Optimization of the expected posterior covariance trace naturally pushes the gliders towards these values that provide the highest information gain.

\subsection{Flow about an unstable node}

\begin{figure}[h!]
\begin{subfigure}{0.49\linewidth}
\begin{center}
\includegraphics[scale=0.35]{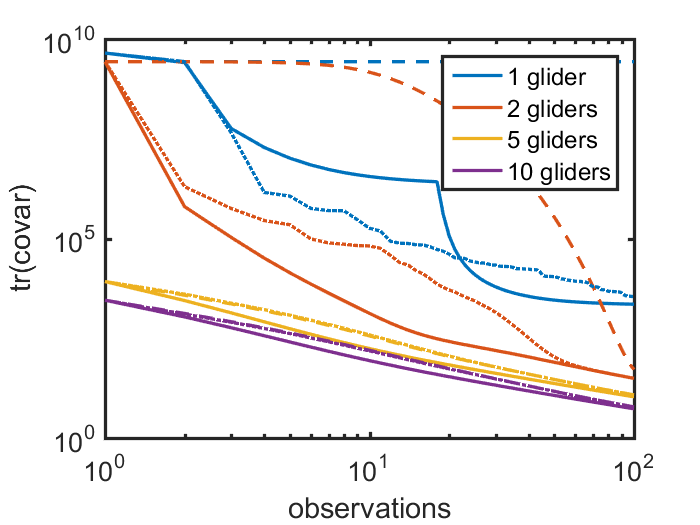}
\caption{\label{f:traceUnstNode:A}Covariance trace, flow about unstable node.}
\end{center}
\end{subfigure}
\begin{subfigure}{0.49\linewidth}
\begin{center}
\includegraphics[scale=0.35]{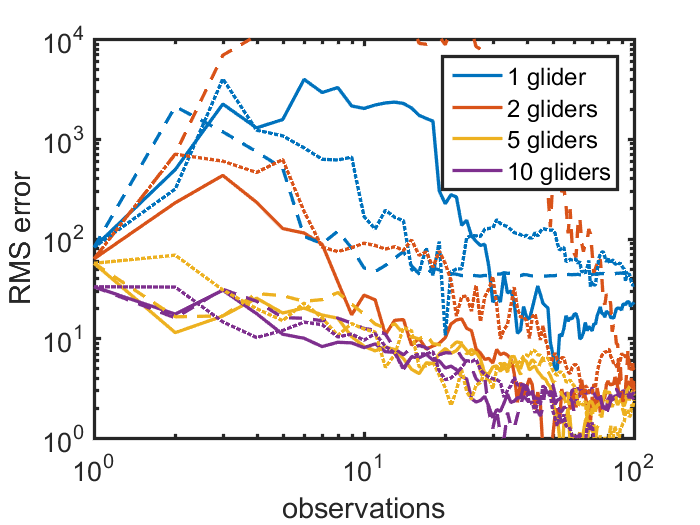}
\caption{\label{f:traceUnstNode:B}RMS error, flow about unstable node.}
\end{center}
\end{subfigure}
\caption{\label{f:traceUnstNode}(A) Decrease in covariance trace with number of observations for flow about unstable node.  For each glider cohort (blue, red, yellow, and purple curves respectively represent 1-, 2-, 5-, and 10-glider cohorts), the solid, dashed, and dotted curves respectively represent locally optimal control, no control, and random control.  Optimal control performs systematically better for sufficiently large number of observations. (B) Same as in (A), but with RMS error used as the performance metric.  Initial growth in the RMS error for the 1- and 2-glider cohorts reflects the high impact of noisy observations on these underresolved cases.}
\end{figure}
\begin{figure}[h!]
\begin{subfigure}{0.49\linewidth}
\begin{center}
\includegraphics[scale=0.35]{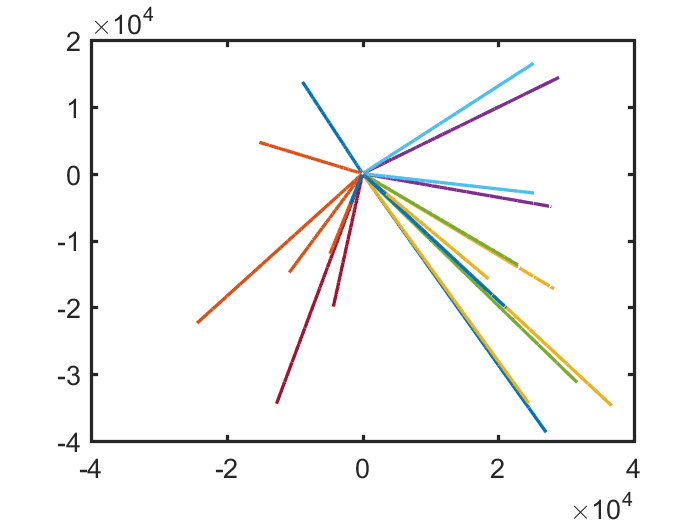}
\caption{\label{f:gliderPathsUnstNode:A}Optimal paths, flow about unstable node.}
\end{center}
\end{subfigure}
\begin{subfigure}{0.49\linewidth}
\begin{center}
\includegraphics[scale=0.35]{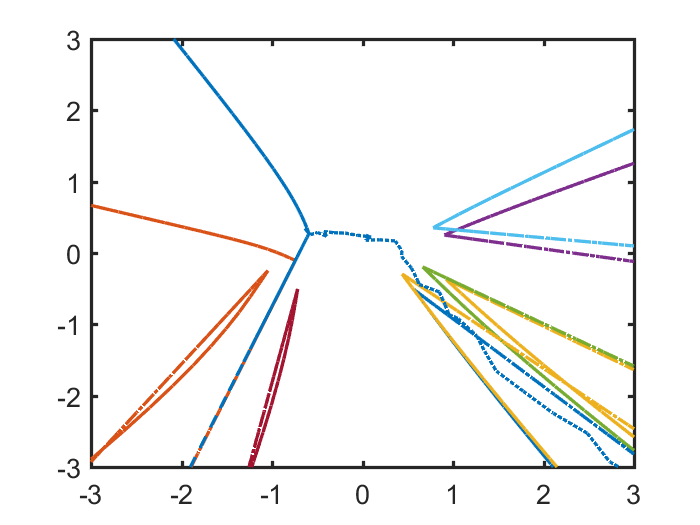}
\caption{\label{f:gliderPathsUnstNode:B}Blowup of (A).}
\end{center}
\end{subfigure}
\caption{\label{f:gliderPathsUnstNode}(A) Optimally controlled glider velocities for flow about unstable node.  Solid lines are optimally controlled, dashed lines have no control, and dotted lines have random control. (B) Same as in (A) but with axes rescaled to give better picture of uncontrolled and randomly controlled paths.}
\end{figure}
The simulations performed around an unstable node use a velocity field given by
\be
v_0 = \begin{pmatrix}1/2\\-1/2\end{pmatrix},\qquad A = \begin{pmatrix}1 & 0\\0 & 1\end{pmatrix},
\ee
such that the center is located at $(-1/2,\,1/2)$.
Figure~\ref{f:traceUnstNode} again shows the advantage of using optimal control, which in this case shows most clearly at intermediate stages of the assimilation process, after sufficient information has been acquired about the flow to produce effective control, and before sufficient observations have been taken to render all methods equally effective.  The exception to this is the case of a single glider, where uncontrolled observations do not successfully learn the flow within the 100 observations taken.  The paths taken by these gliders are shown in Fig.~\ref{f:gliderPathsUnstNode}, where the most salient point is that the controlled gliders follow trajectories that asymptote toward radial lines emanating from the origin for the reason discussed above that this optimizes information gain from observation operator~\ref{e:obsOpLinear}.  The uncontrolled and randomly controlled paths asymptote toward radial lines emanating from the fixed point $(-1/2,\,1/2)$, with the difference between these asymptotic behaviors most noticeable in Fig.~\ref{f:gliderPathsUnstNode:B}.

\subsection{Flow about a saddle}

\begin{figure}[h!]
\begin{subfigure}{0.49\linewidth}
\begin{center}
\includegraphics[scale=0.35]{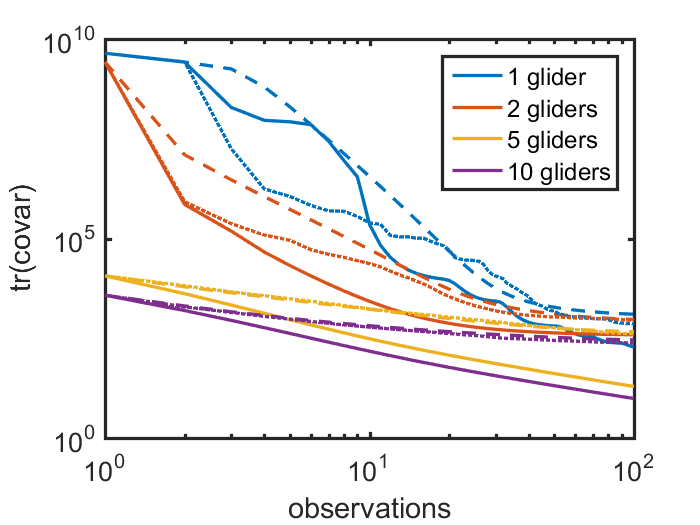}
\caption{\label{f:traceSaddle:A}Covariance trace, flow about saddle.}
\end{center}
\end{subfigure}
\begin{subfigure}{0.49\linewidth}
\begin{center}
\includegraphics[scale=0.35]{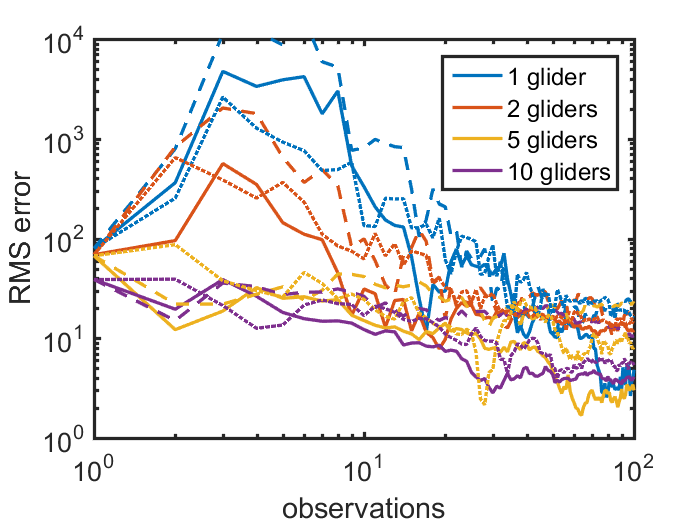}
\caption{\label{f:traceSaddle:B}RMS error, flow about saddle.}
\end{center}
\end{subfigure}
\caption{\label{f:traceSaddle}(A) Decrease in covariance trace with number of observations for flow about saddle.  For each glider cohort (blue, red, yellow, and purple curves respectively represent 1-, 2-, 5-, and 10-glider cohorts), the solid, dashed, and dotted curves respectively represent locally optimal control, no control, and random control.  Optimal control performs systematically better. (B) Same as in (A), but with RMS error used as the performance metric.  Initial growth in the RMS error for the 1- and 2-glider cohorts reflects the high impact of noisy observations on these underresolved cases.}
\end{figure}
\begin{figure}[h!]
\begin{subfigure}{0.49\linewidth}
\begin{center}
\includegraphics[scale=0.35]{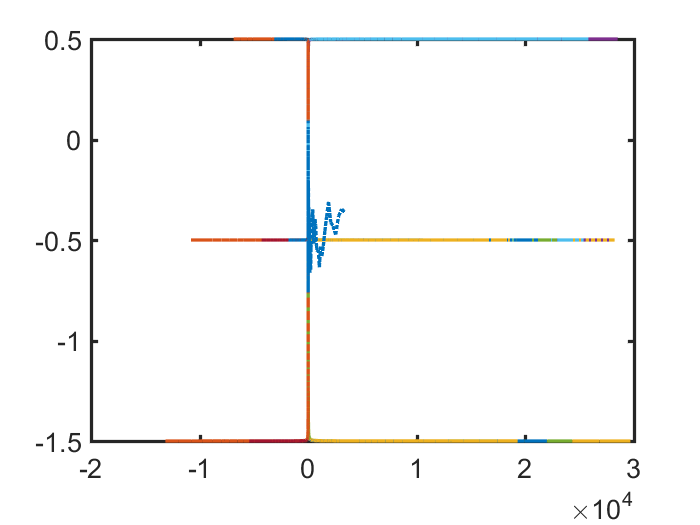}
\caption{\label{f:gliderPathsSaddle:A}Optimal paths, flow about saddle.}
\end{center}
\end{subfigure}
\begin{subfigure}{0.49\linewidth}
\begin{center}
\includegraphics[scale=0.35]{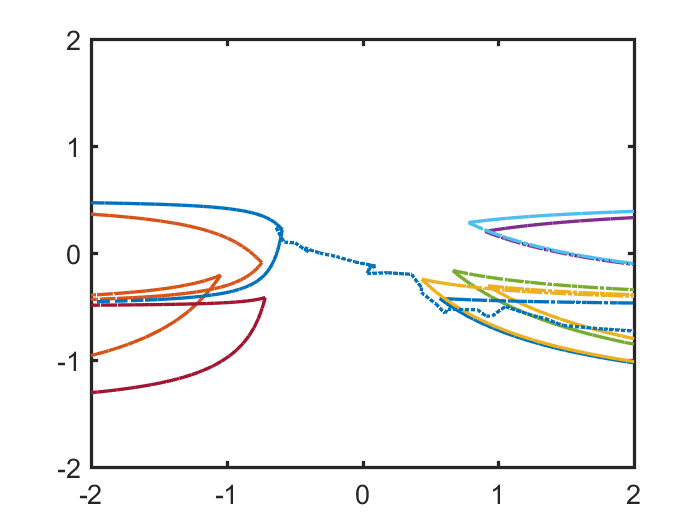}
\caption{\label{f:gliderPathsSaddle:B}Blowup of (A).}
\end{center}
\end{subfigure}
\caption{\label{f:gliderPathsSaddle}(A) Optimally controlled glider velocities for flow about saddle.  Solid lines are optimally controlled, dashed lines have no control, and dotted lines have random control. (B) Same as in (A) but with axes rescaled to give better picture of uncontrolled and randomly controlled paths.}
\end{figure}
The simulations performed around a saddle use a velocity field given by
\be
v_0 = \begin{pmatrix}1/2\\-1/2\end{pmatrix},\qquad A = \begin{pmatrix}1 & 0\\0 & -1\end{pmatrix},
\ee
such that the center is located at $(-1/2,\,-1/2)$.
Figure~\ref{f:traceSaddle} is similar to the cases before, showing systematic improvement in both covariance trace and RMS error when optimal control is used to guide the gliders.  Their paths, shown in Fig.~\ref{f:gliderPathsSaddle}, demonstrate that the optimal control orients itself transverse to the stable manifold in order to maximize the magnitude of the glider's $y$-value.  The uncontrolled and randomly controlled gliders quickly collapse onto the unstable manifold $y=-1/2$.

\subsection{Flow about a stable node}

\begin{figure}[h!]
\begin{subfigure}{0.49\linewidth}
\begin{center}
\includegraphics[scale=0.35]{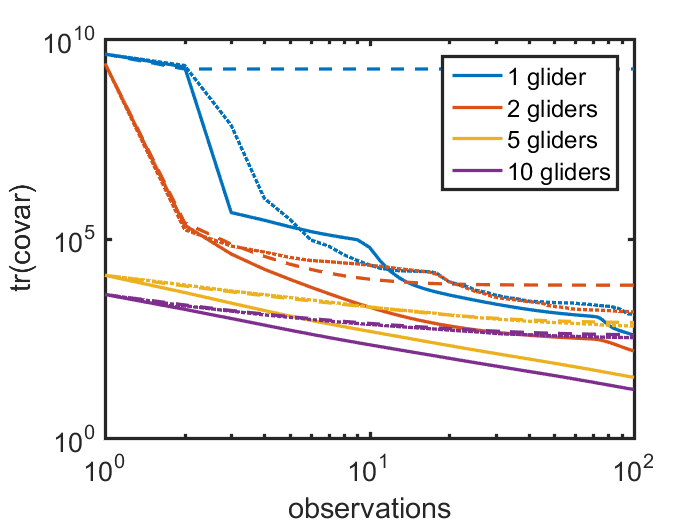}
\caption{\label{f:traceStNode:A}Covariance trace, flow about stable node.}
\end{center}
\end{subfigure}
\begin{subfigure}{0.49\linewidth}
\begin{center}
\includegraphics[scale=0.35]{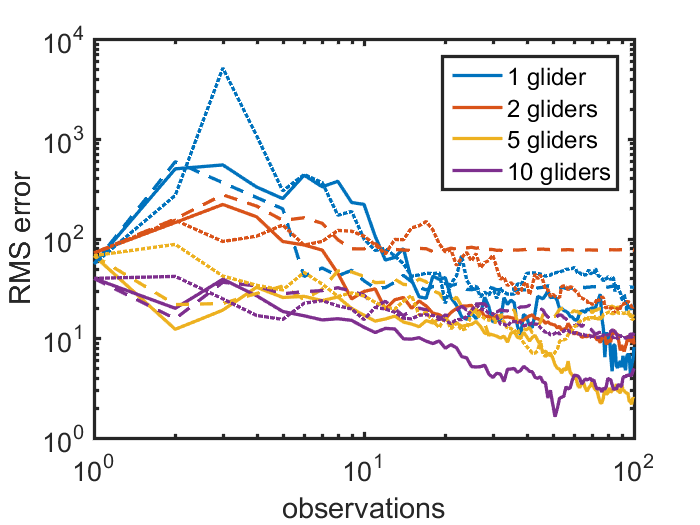}
\caption{\label{f:traceStNode:B}RMS error, flow about stable node.}
\end{center}
\end{subfigure}
\caption{\label{f:traceStNode}(A) Decrease in covariance trace with number of observations for flow about stable node.  For each glider cohort (blue, red, yellow, and purple curves respectively represent 1-, 2-, 5-, and 10-glider cohorts), the solid, dashed, and dotted curves respectively represent locally optimal control, no control, and random control.  Optimal control performs systematically better. (B) Same as in (A), but with RMS error used as the performance metric.  Initial growth in the RMS error for the 1- and 2-glider cohorts reflects the high impact of noisy observations on these underresolved cases.}
\end{figure}
\begin{figure}[h!]
\begin{center}
\includegraphics[scale=0.35]{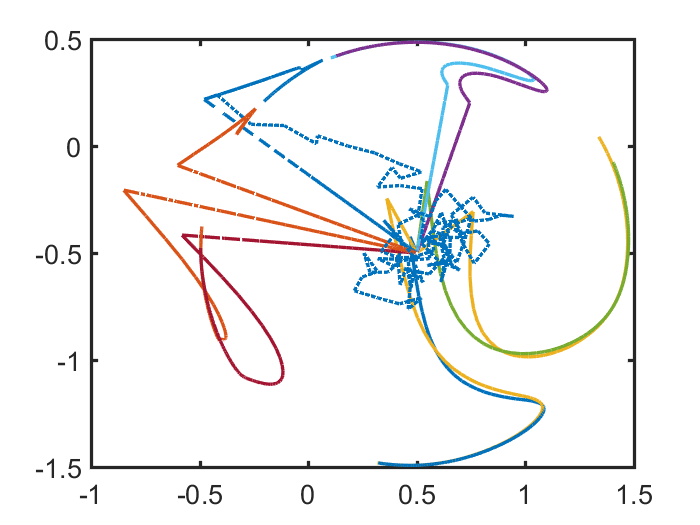}
\caption{\label{f:gliderPathsStNode:A}Optimal paths, flow about stable node.}
\end{center}
\caption{\label{f:gliderPathsStNode}(A) Optimally controlled glider velocities for flow about saddle.  Solid lines are optimally controlled, dashed lines have no control, and dotted lines have random control. (B) Same as in (A) but with axes rescaled to give better picture of uncontrolled and randomly controlled paths.}
\end{figure}
The simulations performed around a stable node use a velocity field given by
\be
v_0 = \begin{pmatrix}1/2\\-1/2\end{pmatrix},\qquad A = \begin{pmatrix}-1 & 0\\0 & -1\end{pmatrix},
\ee
such that the center is located at $(1/2,\,-1/2)$.
In this case, Figure~\ref{f:traceStNode} demonstrates how, in the absence of control, observations from the gliders are eventually rendered uninformative as the gliders collapse onto the fixed point.  Random control offers some improvement, but optimal control is systematically better, particularly once sufficient information about the flow has been acquired.  The paths shown in Fig.~\ref{f:gliderPathsStNode} demonstrate that the optimal control strives to keep the gliders from collapsing into the attractor.

\section{Discussion}
\label{s:discussion}

The methodology presented above provides a systematic way to iteratively improve the estimate of a velocity field obtained through the use of observations made by controlled autonomous vehicles.  We have shown through the use of simple examples consisting of stationary linear flows near hyperbolic fixed points with different stability properties that this approach is effective, providing significant improvement over uncontrolled AVs and over AVs that are driven at maximal speed in random directions between each assimilation stage.  We obtain the locally optimal paths through solution of the Euler-Lagrange equations resulting from an objective function provided by the expected variance of the posterior distribution which, in the case of direct velocity observations, depends on the locations of the AVs at each assimilation time.  As is often the case in optimal control problems, these equations are not always solvable, and we provide a regularization method that provably works in the case of trivial flows and has resolved all cases of nonexistence of solutions in the linear flows considered here.

The example flows considered here are simple linear flows for the purpose of clear illustration of the optimization results.  We now consider the natural extensions to this method towards application to more realistic flows, all of which are the subject of ongoing work.

\subsection{Nonlinear, random, and time-dependent flows}

\begin{figure}[h!]
\begin{subfigure}{0.49\linewidth}
\begin{center}
\includegraphics[scale=0.4]{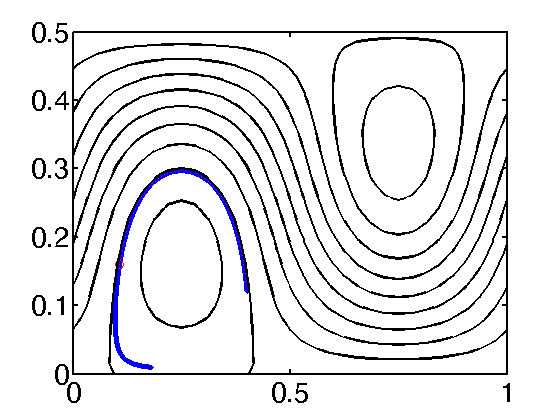}
\caption{\label{f:pathVariance:A}Uncontrolled path, double gyre.}
\end{center}
\end{subfigure}
\begin{subfigure}{0.49\linewidth}
\begin{center}
\includegraphics[scale=0.4]{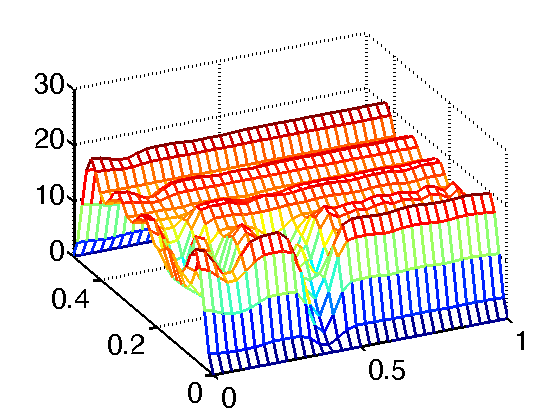}
\caption{\label{f:pathVariance:B}Variance, uncontrolled path in double gyre.}
\end{center}
\end{subfigure}
\begin{subfigure}{0.49\linewidth}
\begin{center}
\includegraphics[scale=0.4]{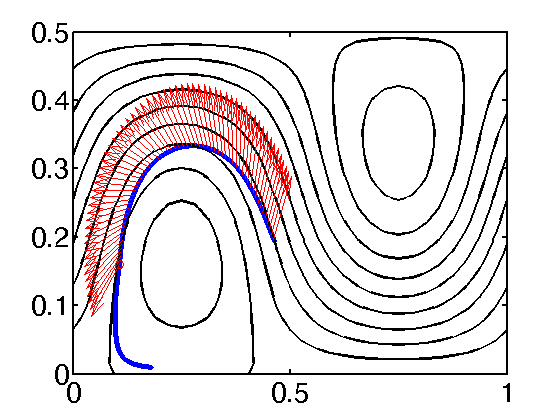}
\caption{\label{f:pathVariance:C}Controlled path, double gyre.}
\end{center}
\end{subfigure}
\begin{subfigure}{0.49\linewidth}
\begin{center}
\includegraphics[scale=0.4]{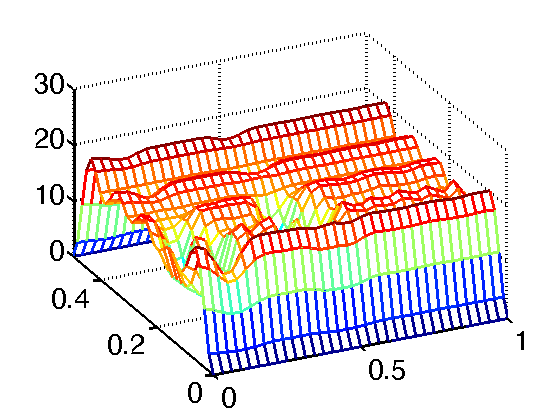}
\caption{\label{f:pathVariance:D}Variance, controlled path in double gyre.}
\end{center}
\end{subfigure}
\caption{\label{f:pathVariance}(A) Typical uncontrolled path (blue) trapped in gyre over two evolution periods with an observation (red circle) between them.  Black curves are streamlines of true flow.  (B) Variance in estimate $\tilde\psi(x,y)$ after two observations along uncontrolled path from (A). (C) Typical controlled path (blue), with control vector depicted in red. (D) Variance in $\tilde\psi(x,y)$ after two observations along controlled path from (C).}
\end{figure}

A primary motivation for the work described in~\cite{mcdougall_jones_2015} was the observation that Lagrangian {\em drifters\/} (i.e., with no capacity for self-propulsion) often encounter flow barriers that prevent their effectiveness in learning regional flows.  The kinematic, incompressible flow~\cite{samelson_lagrangian_2006} used as a test case in that work took the form of a stream function given by
\begin{equation}
  \psi = -\pi y + \sin(2 \pi x) \sin(2 \pi y),
  \quad (x, y) \in [0, 1] \times [0, 0.5],
\end{equation}
with associated (stationary) velocity field
\begin{align}
  v_1 &= -\frac{\partial \psi}{\partial y} = \pi - 2 \pi \sin(2 \pi x) \cos(2 \pi y) \\
  v_2 &=  \frac{\partial \psi}{\partial x} =       2 \pi \cos(2 \pi x) \sin(2 \pi y).
\end{align}
If the mean flow $\langle v \rangle =(\pi,0)^T$ is assumed known with the rest assumed to satisfy periodic boundary conditions in $x$ and Dirichlet boundary conditions in $y$, an estimate $\tilde{\psi}$  of the stream function can be represented through a finite Fourier basis given by
\be
\tilde{\psi} = \pi y + \sum_{n=1}^N a_{on}\sin(2n\pi y) +
\sum_{m=1}^M\sum_{n=1}^N[a_{mn}\cos(2m\pi x)+b_{mn}\sin(2m\pi x)]\sin(2n\pi y).
\ee
Expressing the state vector $\hat{\psi}$ as a column of these coefficients with dimension $N(1+2M)$, the observation operator mapping from state space to observation space then takes the form
\be
H = \begin{pmatrix}
k_1\cos(k_1y_1) & \dots & k_1\cos(\mu_1x_1)\cos(k_1y_1) & k_1\sin(\mu_1x_1)\cos(k_1y_1) & \dots\\
0 & \dots & \mu_1\sin(\mu_1x_1)\sin(k_1y_1) & -\mu_1\cos(\mu_1x_1)\sin(k_1y_1) & \dots\\
k_1\cos(k_1y_2) & \dots & & & \\
\vdots & & & &\end{pmatrix},
\ee
from which one can see the dependence on glider positions $z_k=(x_k,y_k)^T$.  More complex flows imply more significant difficulties with nonexistence of solutions to the Euler-Lagrange equations, but preliminary results are nevertheless promising.  Depicted in Fig.~\ref{f:pathVariance} are comparisons of controlled and uncontrolled paths, with their related reduction in the spatial variance of the stream function estimate.  Whereas the controlled path is trapped within the gyre, producing observations that do little to reduce uncertainty outside of the gyre, the optimal control algorithm produces a path that crosses the separatrix into the meandering flow, producing a corresponding reduction in uncertainty there.  This work is ongoing, as is the incorporation of weakly nonlinear time evolution models for the flow that fit naturally into the framework of the extended Kalman filter.

\subsection{Observations of position}

As described in Sec.~\ref{s:setup}, direct observations of the velocity local to a glider produce an observation operator of the form given in Eqn.~\ref{e:obsOpVel}, where it is clear how the glider positions at assimilation time influence the Kalman gain and therefore the expected posterior covariance.  When it is the gliders' positions that are observed, the methodology proposed in Ref.~(\cite{Kuznetsov2003}) poses as the state the augmented vector consisting of a parametric representation of the velocity field concatenated with the glider positions themselves.  This produces the observation operator given in Eqn.~\ref{e:obsOpPos}, where it is no longer immediately clear how to formulate the dependence of expected posterior variance on glider path in such a way as to compute an optimal control.  In this case, since the estimated state includes the glider positions themselves, the covariance of the state must be evolved between assimilation times, and this time evolution of the covariance is where the control plays a role.  The extension of our control algorithm to this case will be presented elsewhere.

\bibliographystyle{elsarticle-num}
\bibliography{paper}

\end{document}